\documentclass[12pt]{amsart}

\usepackage{amssymb}
\usepackage{amsmath}
\usepackage{latexsym}
\usepackage{epsfig}
\usepackage{graphics}
\usepackage{color}

\newcommand{\N}{\ensuremath{\mathbb{N}}}
\newcommand{\A}{\ensuremath{\mathbb{A}}}
\newcommand{\Q}{\ensuremath{\mathbb{Q}}}
\newcommand{\Z}{\ensuremath{\mathbb{Z}}}

\newcommand{\F}{\ensuremath{\mathbb{F}}}

\newcommand{\Zp}{\ensuremath{{\mathbb Z}_p}}

\DeclareMathOperator{\ord}{ord}

\DeclareMathOperator{\Frac}{Frac}
\DeclareMathOperator{\ac}{ac}

\DeclareMathOperator{\smod}{mod}

\DeclareMathOperator{\Spec}{Spec}
\DeclareMathOperator{\Min}{Min}

\newcommand\cO{\mathcal{O}}

\newcommand\fX{\mathfrak{X}}

\newcommand\frm{\mathfrak{m}}

\newtheorem{mainTheorem}[subsection]{Main Theorem}
\newtheorem{lemma}[subsection]{Lemma}
\newtheorem{definition}[subsection]{Definition}  

\newtheorem{remark}[subsection]{Remark}

\newtheorem{sectionItem}[subsection]{}

\newtheorem{namedTheorem}[subsection]{\theoremname}
\newcommand{\theoremname}{testing}
\newenvironment{named}[1]{\renewcommand\theoremname{#1}
\begin{namedTheorem}}
{\end{namedTheorem}}

\begin{document}

\title{PROOF OF A CONJECTURE OF COLLIOT-TH\'{E}L\`{E}NE}

\author{Jan Denef}

\address{KULeuven, Department of Mathematics, Celestijnenlaan 200 B,
B-3001 Leuven (Heverlee), Belgium.}
\email{jan.denef@wis.kuleuven.be}


\date{January 22, 2016}

\begin{abstract}
We prove a conjecture of Colliot-Th\'{e}l\`{e}ne that implies the Ax-Kochen Theorem on $p$-adic forms. We obtain it as an easy consequence of a diophantine excision theorem whose proof forms the body of the present paper.
\end{abstract}

\maketitle

\section{Introduction}

\noindent
In this paper we prove the following conjecture of Colliot-Th\'{e}l\`{e}ne \cite{Colliot-1}.

\begin{sectionItem} {\rm\textbf{Colliot-Th\'{e}l\`{e}ne's Conjecture.}} \label{conj-Colliot}
Let $f: X \rightarrow Y$ be a dominant morphism of nonsingular proper geometrically integral varieties over $\Q$, with geometrically integral generic fibre. Assume that for any nontrivial discrete valuation on the function field K of $Y$,
with valuation ring $A \supset \Q$, there exists an integral regular $A$-scheme $\fX$, flat and proper over $A$, with generic fibre $K$-isomorphic to the generic fibre of $f$, and special fibre having an irreducible component of multiplicity 1 which is geometrically integral. Then the map $X(\Q_p) \rightarrow Y(\Q_p)$, induced by f, is surjective for almost all primes $p$.
\end{sectionItem}

\noindent
Here $\Q_p$ denotes the field of $p$-adic numbers, and with \textquotedblleft almost all primes\textquotedblright $\,$ we mean \textquotedblleft  all but a finite number of primes\textquotedblright.

Actually we prove a stronger result, namely:
\begin{mainTheorem}\label{th-main}
Let $f: X \rightarrow Y$ be a dominant morphism of nonsingular proper geometrically integral varieties over $\Q$, with geometrically integral generic fibre. Assume for any modification $f': X' \rightarrow Y'$ of $f$, with the same generic fibre as $f$, and $X', Y'$ nonsingular, and for any prime divisor $D'$ on $Y'$, the following: the divisor $f'^{*}(D')$ on $X'$ has an irreducible component $C'$ with multiplicity 1 and geometrically integral generic fibre over $D'$ (i.e. the morphism $C' \rightarrow D'$, induced by $f'$, has geometrically integral generic fibre). Then the map $X(\Q_p) \rightarrow Y(\Q_p)$, induced by f, is surjective for almost all primes $p$.
\end{mainTheorem}

We say that $f'$ is a \emph{modification} of $f$ if $f'$ fits into a commutative square of morphisms of varieties, with vertical arrows $f, f'$, and horizontal arrows birational proper morphisms $X' \rightarrow X, \; \; Y' \rightarrow Y$, see Definition \ref{def-modif}.

The conjecture of Colliot-Th\'{e}l\`{e}ne is a direct consequence of Theorem \ref{th-main}, because any $D'$ as in the theorem induces a discrete valuation on the function field of $Y'$, which equals the function field of $Y$. Such a discrete valuation on the function field of $Y$ is called a \emph{divisorial valuation}. Moreover, if for this valuation there exists an $A$-scheme $\fX$ as in the conjecture, then the special fibre of any integral regular proper flat $A$-scheme $\fX'$, with the same generic fibre as $\fX$, has an irreducible component of multiplicity 1 which is geometrically integral. Indeed this is Proposition 3.9.(b) in Colliot-Th\'{e}l\`{e}ne's lecture notes \cite{Colliot-2}. Thus the hypotheses in the statement of the conjecture of Colliot-Th\'{e}l\`{e}ne imply the hypotheses in the statement of Theorem \ref{th-main}, even if we restrict the assumption in the conjecture to divisorial valuations.

Note that Theorem \ref{th-main} is substantially stronger than the conjecture of Colliot-Th\'{e}l\`{e}ne, because it implies that the assumption in the conjecture is only required for divisorial discrete valuations on the function field of $Y$.

Colliot-Th\'{e}l\`{e}ne \cite{Colliot-1} proved the following: if $f: X \rightarrow Y$ is the universal family over $\Q$ of all projective hypersurfaces of degree $d$ in projective $n$-space, with $n \geq d^{\, 2}$, then $f$ satisfies the hypotheses of the conjecture and hence also the hypotheses of Theorem \ref{th-main}. Since our proof of Theorem \ref{th-main} is purely algebraic geometric, this yields a new proof of the Theorem of Ax and Kochen \cite{AxKoch-1} on p-adic forms, that does not rely on methods from mathematical logic. The Theorem of Ax and Kochen states that for each $d\in \N$ there exists $N \in \N$ such that for all primes $p > N$, each hypersurface of degree $d$ in projective $n$-space over $\Q_p$, with $n \geq d^{\, 2}$, has a $\Q_p$-rational point. 

One of the motivations of Colliot-Th\'{e}l\`{e}ne in formulating his conjecture was to obtain an algebraic geometric proof of the Ax-Kochen Theorem that, unlike all previous ones, does not rely on methods from mathematical logic. At the same time, the author of the present paper also found another purely algebraic geometric proof of the Ax-Kochen Theorem, see \cite{Denef-AxKochen}. Both proofs are based on the Tameness Theorem (see section \ref{sec-tameness}), which is proved in \cite{Denef-AxKochen} using the Weak Toroidalization Theorem of Abramovich an Karu \cite{Abramovich-Karu} (extended to non-closed fields \cite{Abramovich-Denef-Karu}).

We prove the Main Theorem \ref{th-main} in section \ref{sec-proofMainThm}, as an easy consequence of what we call a diophantine excision theorem . The proof of this Diophantine Excision Theorem \ref{th-purity} forms the body of the present paper and is contained in section \ref{sec-purity}. It depends on the Tameness Theorem \ref{th-tameness}, which is treated in section \ref{sec-tameness}. Using mathematical logic one can give a simpler proof of Colliot-Th\'{e}l\`{e}ne's Conjecture \ref{conj-Colliot}. However we don't see how to extend this to prove the stronger Theorem \ref{th-main} or the Diophantine Excision Theorem \ref{th-purity}. This alternative proof is given in section \ref{prf-alternative}.
Preliminaries about modifications of morphisms and multiplicative residues are given in sections \ref{sec-modif} and \ref{SectionMultRes}.

A previous version of the present paper was posted on the ArXiv in 2011. In that version the Diophantine Excision Theorem was called the Diophantine Purity Theorem.

\begin{sectionItem} Terminology and notation. \label{conventions} \rm

For any prime $p$ we denote the ring of $p$-adic integers by $\Z_p$, the field of $p$-adic numbers by $\Q_p$, and the field with $p$ elements by $\F_p$. The $p$-adic valuation on $\Q_p$ is denoted by $\ord_p$.

In the present paper, $R$ will always denote a noetherian integral domain. A \emph{variety over} $R$ is an integral separated scheme of finite type over $R$.  With a \emph{morphism} of varieties over $R$ we mean an $R$-morphism of schemes over $R$.
A rational function $x$ on a variety $X$ over $R$ is called \emph{regular at} a point $P \in X$ if it belongs to the local ring $\cO_{X,P}$ of $X$ at $P$, and it is called \emph{regular} if it is regular at each point of $X$.

\emph{Uniformizing parameters over} $R$ on a variety $X$ over $R$, are regular rational functions on $X$ that induce an \'{e}tale morphism from $X$ to an affine space over $R$.

A \emph{reduced strict normal crossings divisor over} $R$ on a smooth variety $X$ over $R$ is a closed subset $D$ of $X$ such that for any $P \in X$ there exist uniformizing parameters $x_1, \dots, x_n$ over $R$ on an open neighborhood of $P$, such that for any irreducible component $C$ of $D$, containing $P$, there is an $i \in \{1, \dots, n\}$ which generates the ideal of $C$ in $\cO_{X,P}$.

\end{sectionItem}

For ease of notation we work with the completions $\Q_p$ of $\Q$, but all results in the present paper remain true replacing $\Q$ by any number field $K$, and the completions $\Q_p$ by the non-archimedean completions of $K$. \\

\noindent
\emph{Acknowledgments.} We thank J.-L. Colliot-Th\'{e}l\`{e}ne, O. Gabber, and O. Wittenberg for a mutual conversation that resulted in the alternative proof \ref{prf-alternative}. We also thank D. Abramovich, R. Cluckers, J.-L. Colliot-Th\'{e}l\`{e}ne, S. D. Cutkosky, K. Karu, and O. Wittenberg for stimulating conversations and useful information. Moreover Cao's paper \cite{Cao-ColliotConjecture} was helpful while producing the final version of the present paper.

\section{Modifications of Morphisms} \label{sec-modif}
\begin{definition} \label{def-modif} \rm Let $R$ be a noetherian integral domain, and $X$ a variety over $R$. A \emph{modification of} $X$ is a proper birational morphism $X' \rightarrow X$ of varieties over $R$.

Let $f: X \rightarrow Y$ be a dominant morphism of varieties over $R$. A \emph{modification of} $f$ is a morphism $f': X' \rightarrow Y'$ of varieties over $R$, which fits into a commutative diagram
$$X' \,\xrightarrow{\alpha }\, X $$
$$f' \downarrow \quad \;\; \; \; \downarrow f$$
$$Y' \,\xrightarrow{\beta }\, Y $$
with $\alpha$ a modification of $X$, and $\beta$ a modification of $Y$. This implies that $f'$ is dominant. Clearly, if $f$ is proper, then also $f'$ is proper.

When $f: X \rightarrow Y$ is a dominant morphism of varieties over $R$, and $\beta: Y' \rightarrow Y$ is a modification of $Y$, then there exists a unique irreducible component $X'$ of the fibre product $Y' \times_Y X$ that dominates $X$. Let $f'$ and $\alpha$ be the restrictions to $X'$ of the projections $Y' \times_Y X \rightarrow Y'$, and $\, Y' \times_Y X \rightarrow X$. Then $\alpha$ is a modification of $X$, and we call $f'$ the \emph{strict transform of $f$ with respect to} $\beta$. Clearly $f'$ is a modification of $f$. Such a modification is called a \emph{strict modification of} $f$. Note that any strict modification of $f'$ is also a strict modification of $f$.
\end{definition}

\begin{sectionItem}\rm {\textbf{Observations}} \label{obs-relatedToModifications}

(a) Let $f: X \rightarrow Y$ be a morphism of schemes of finite type over an excellent henselian discrete valuation ring $R$, with $Y \otimes_R \Frac(R)$ nonsingular, where $\Frac(R)$ denotes the fraction field of $R$. Let $S$ be a closed subscheme of $Y$, containing no irreducible component of $Y$.
If $\,Y(R) \setminus S(R)\, \subset \, f(X(R))$, then $Y(R) \subset f(X(R))$. Indeed this follows from Greenberg's theorem \cite{Greenberg}, because $Y(R) \setminus S(R)$ is dense in $Y(R)$ with respect to the adic topology on $Y(R)$, since $Y \otimes_R \Frac(R)$ is nonsingular.

(b) Let $f: X \rightarrow Y$ be a dominant morphism of varieties over $\Z$, and let $f': X' \rightarrow Y'$ be a strict modification of $f$. Assume that $Y\otimes \Q$ and $Y'\otimes \Q$ are nonsingular, and let $p$ be a prime number. Then, the map $X(\Zp) \rightarrow Y(\Zp)$, induced by $f$, is surjective, if and only if the map $X'(\Zp) \rightarrow Y'(\Zp)$, induced by $f'$, is surjective. This remains true when $f'$ is a modification of $f$ which is not strict, if we assume that also $X\otimes \Q$ is nonsingular. These claims follow directly from (a). Indeed, by (a) and the valuative criterium for properness, any modification of a variety $V$ over $\Z$, with $V \otimes \Q$ nonsingular, induces a surjection on $\Z_p$-rational points.

\end{sectionItem}

\noindent
\emph{Remark.} We will often use (without mentioning) the following well known facts. Any morphism $f_0: X_0 \rightarrow Y_0$ of varieties over $\Q$ has a model $f$ over $\Z$.
This means that $f$ is a morphism $f: X \rightarrow Y$ of varieties over $\Z$ whose base change to $\Q$ is isomorphic to $f_0$. Combining this with Nagata's compactification theorem (see e.g. \cite{Lutke-Compactification}), we see that we can choose $f$ to be proper, when $f_0$ is proper. Two models of $f_0$ over $\Z$ become isomorphic after base change to $\Z[1/N]$, for some positive integer $N$. Hence, if $f_0$ is proper and $f$ is a model of $f_0$ over $\Z$, then $f \otimes \Z[1/N]$ is proper for some $N \in \N$.

\section {Multiplicative residues} \label{SectionMultRes}

Let $R$ be a noetherian integral domain, and $X$ a variety over $R$.
Let $A$ be any local $R$-algebra which is an integral domain. We denote by $\frm_A$ its maximal ideal, by $\Frac(A)$ its field of fractions, and by $\eta_A$ the generic point of $\Spec(A)$. For any $A$-rational point $a  \in X(A)$ on $X$ we denote by $a \!\mod \frm_A$ the $A / \frm_A$-rational point on $X$ induced by $a$. For any $x \in \cO_{X,a(\eta_A)}$ the pullback $a^*(x)$ of $x$ to $\Frac(A)$ is denoted by $x(a) \in \Frac(A)$. Moreover, for $a, a' \in X(A)$ we write $a \equiv a' \!\mod \frm_A$ to say that  $a \!\mod \frm_A \;=\; a' \!\mod \frm_A$.

\begin{named}{Definition} \label{DefinitionMultRes1} \rm
Let $z,z' \in \Frac(A)$. The elements $z,z'$ have the \emph{ same multiplicative residue} if
$$ z' \in z(1 + \frm_A).$$
\newline
Let $a,a' \in X(A)$ and let $x_1, \dots , x_r$ be rational functions on $X$. The points $a,a'$ have the \emph{same residues
with respect to} $x_1, \dots , x_r$ if $a \equiv a'\!\mod\frm_A$ and, for $i = 1, \dots, r$, the following two conditions hold.
\begin{enumerate}
\item The rational function $x_i$ is regular at $a(\eta_A)$ if and only if it is regular at  $a'(\eta_A)$.
\item $x_i(a), x_i(a') \in \Frac(A)$ have the same  multiplicative residue if $x_i$ is regular at both $a(\eta_A)$ and $a'(\eta_A)$.
\end{enumerate}
\end{named}

The following lemma also appears in \cite{Denef-AxKochen}.

\begin{lemma}
\label{LemmaResidues1}
Let $X$ be an affine variety over $R$, and let $x_1, \dots, x_r$ be rational functions on $X$. Then there exist regular rational functions $x'_1, \dots, x'_{r'}$ on $X$ such that for any local R-algebra $A$, which is an integral domain, and any $a, a' \in X(A)$ we have the following. The points $a$ and $a'$ have the same residues with respect to $x_1, \dots, x_r$ if they have the same residues with respect to  $x'_1, \dots, x'_{r'}$.
\end{lemma}

\emph{Proof}. This is clear, by taking for $x'_1, \dots, x'_{r'}$ any finite list of regular rational functions on $X$ which satisfies the following condition. For each $i \in \{1, \dots, r\}$ and each $P \in X$ with $x_i$ regular at $P$, there are elements $x'_j$ and $x'_k$ in this list with $x_i = x'_j/x'_k$, and $x'_k(P)\neq 0$.  Obviously, such a finite list exists if $X$ is affine. $\square$

\begin{lemma}
\label{LemmaResidues2}
Let $X$ be a variety over $R$, and let $z, x_1, \dots, x_r$ be regular rational functions on $X$. Assume that $z$ can be written as a unit in $\Gamma(X, \mathcal{O}_X)$ times a monomial in the $x_i$. Then for any local R-algebra $A$, which is an integral domain, and any $a, a' \in X(A)$ we have the following. The points $a$ and $a'$ have the same residues with respect to $z$ if they have the same residues with respect to  $x_1, \dots, x_r$.
\end{lemma}

\emph{Proof}. Obvious, and left to the reader. $\square$

\begin{definition} \rm
Let $X$ be a variety over $\Z$, and $x_1,\dots, x_r$ regular rational functions on $X$. Let $z = (z_1,\dots, z_r) \in \Z_p^r$. We say that \emph{the multiplicative residue of $z$ is realizable with respect to} $x_1,\dots, x_r$ if there exists $a \in X(\Z_p)$ such that $x_i(a)$ and $z_i$ have the same multiplicative residue for each $i = 1, \dots, r$.
\end{definition}

\begin{definition} \rm
For any $w \in \Q_p$, the \emph{angular component modulo $p$ of $w$} is defined as
$$ \overline{\ac}(w) := w p^{-\ord_p(w)}\, \smod p \, \in \F_p \, ,$$
with the convention that $\overline{\ac}(0) := 0$.

Note that any $w, w' \in \Q_p$ have the same multiplicative residue if and only if they have the same $p$-adic valuation and the same angular component modulo $p$.
\end{definition}

The following rather technical lemma will be used in the proof of the Surjectivity Criterium \ref{surjectCriterium}. It is a direct consequence of the Theorem of Pas \cite{Pas-1} on uniform $p$-adic quantifier elimination. The work of Pas is based on methods from mathematical logic. Below we give a purely algebraic geometric proof of this lemma which is based on embedded resolution of singularities.

\begin{lemma}\label{lem-forSurjectCrit}
Let $X$ be a variety over $\Z$, and $x_1,\dots, x_r$ regular rational functions on $X$. There exists a finite partition of $\,\N^r$ such that for almost all primes $p$ we have the following. Let $z, z' \in \Z_p^r\,$. Assume that the $p$-adic valuations of $z$ and $z'$ are in a same stratum of the partition, and that $\,\overline{\ac}(z_i) = \overline{\ac}(z'_i)\,$ for each $i$.
Then the multiplicative residue of $z$ is realizable with respect to $x_1,\dots, x_r$, if and only if the same holds for $z'$.
\end{lemma}
\noindent
\emph{Proof.} Let $D$ be the union of the zero loci of the regular rational functions $x_1, \dots, x_r$ on $X$, considered as a subset of $X$. Using embedded resolution of singularities of $D  \otimes \Q \subset X \otimes \Q$, and induction on the dimension of $X \otimes \Q$, modifying $X$ and  inverting a finite number of primes, we may assume the following. The variety $X$ is smooth over $\Z$, and $D$ is a reduced strict normal crossings divisor over $\Z$ (in the sense of subsection \ref{conventions}). This reduction is easily verified applying the valuative criterion of
properness to the resolution morphism and using the induction hypothesis to take care of the exceptional locus of the resolution. Hence, by covering $X$ with finitely many suitable open subschemes, we can further assume that $X$ is affine, and that each $x_i$ can be written as a unit $u_i$ in $\Gamma(X, \mathcal{O}_X)$ times a monomial in uniformizing parameters $y_1, \dots,  y_n$ over $\Z$ on $X$. As recalled in subsection \ref{conventions}, this means that $y_1, \dots,  y_n$ induce an \'{e}tale morphism from $X$ to affine $n$-space over $\Z$. 

Let $E$ be the matrix over $\Z$, with $r$ rows and $n$ columns, consisting of the exponents of these monomials, and let $\Delta$ be the linear map $\Z^n \rightarrow \Z^r$ determined by the matrix $E$.
\newline
For each subset $S \subset \{1, \dots, n \}$, set
\begin{align*}
&\Gamma_{S} := \{(\alpha_1, \dots,\alpha_n) \in \N^n \,|\, \forall j:\, \alpha_j = 0 \Leftrightarrow j \in  S\}.
\end{align*}
Choose a finite partition of $\N^r$ such that each $\Delta(\Gamma_{S})$ is a union of strata.
Let $z, z' \in \Z_p^r\,$ be as in the lemma and assume that the multiplicative residue of $z$ is realizable with respect to $x_1,\dots, x_r$ by an element $a \in X(\Z_p)$. We have to show that $z'$ is also realizable. By slightly moving $a$, we may suppose that $y_j(a) \not = 0$ for all $j=1, \dots, n$. Let $S \subset \{1, \dots, n \}$ be such that $(\ord_p y_1(a), \dots, \ord_p y_n(a)) \in \Gamma_{S}$. Hence $\ord_p z = (\ord_p x_1(a), \dots, \ord_p x_r(a)) \in \Delta(\Gamma_{S})$, since $u_i(a)$ is a unit in $\Z_p$ for each $i$. Because the $p$-adic valuations of $z$ and $z'$ are in a same stratum of the partition, also $\ord_p z'$ is an element of $\Delta(\Gamma_{S})$. Hence there exists $ \alpha' = (\alpha'_1, \dots,\alpha'_n) \in \Gamma_{S}$ with $\Delta(\alpha') = \ord_p(z')$. Note that $\alpha'_j = 0$ if and only if $\ord_p(y_j(a)) = 0$, because $\alpha' \in \Gamma_{S}$. 

By Hensel's Lemma, applied to the \'{e}tale morphism induced by the $y_1, \dots, y_n$, there exists $a' \in X(\Z_p)$ with $\,a'\,\smod p \, = \, a\,\smod p \,$ and $\ord_p (y_j(a')) = \alpha'_j$ and $\overline{\ac}(y_j(a')) = \overline{\ac}(y_j(a))$, for $j = 1, \dots, n$. Indeed any element of $\Z_p^n$ which is congruent mod $p$ to the image of $a$ under this \'{e}tale morphism, can be lifted to a point $a' \in X(\Z_p)$ congruent to $a$. Now we have that $\overline{\ac}(x_i(a')) = \overline{\ac}(x_i(a)) = \overline{\ac}(z_i) = \overline{\ac}(z'_i)$ and
$\ord_p (x_i(a')) = (\Delta(\alpha'))_i = \ord_p (z'_i)$, because $u_i(a)$ and $u_i(a')$ are units in $\Z_p$ which are congruent mop $p$. Hence $x_i(a')$ and $z'_i$ have the same multiplicative residue for $i = 1, \dots, r$. Thus the multiplicative residue of $z'$ with respect to $x_1,\dots, x_r$ is realized by $a'$. This terminates the proof of the lemma. $\square$

\section{Tameness and the Surjectivity Criterium} \label{sec-tameness}
\noindent
The following result is a special case of the Tameness Theorem of \cite{Denef-AxKochen} (together with Remark 5.2 of \cite{Denef-AxKochen}).
\begin{named}{Tameness Theorem} \label{th-tameness}
Let $f : X \rightarrow Y$ be a morphism of varieties over $\Z$.
Given rational functions $x_1, \dots, x_r$ on $X$, there exist rational functions $y_1, \dots, y_s$ on $Y$, such that for almost all primes $p$ we have the following.
Any $b \in Y(\Z_p)$ having the same residues with respect to $\, y_1, \dots, y_s$ as an image $f(a')$, with $a' \in X(\Z_p)$, is itself an image of an $a \in X(\Z_p)$ with the same residues as $a'$ with respect to $\, x_1, \dots, x_r$.
\newline
Moreover, if $Y$ is affine, then we can choose $y_1, \dots, y_s$ to be regular rational functions on $Y$.
\end{named}

\noindent
This special case, and the more general result in \cite{Denef-AxKochen}, can be proved easily by using Basarab's theorem \cite{Basarab-1} on elimination of quantifiers. The special case itself is also an easy consequence of the theorem of Pas \cite{Pas-1} on uniform $p$-adic quantifier elimination. The works of Pas and Basarab are based on methods from mathematical logic. However in \cite{Denef-AxKochen} we gave a purely algebraic geometric proof of the Tameness Theorem which is based on the Weak Toroidalization Theorem of Abramovich an Karu \cite{Abramovich-Karu} (extended to non-closed fields \cite{Abramovich-Denef-Karu}).

We briefly sketch the geometric proof of the Tameness Theorem of \cite{Denef-AxKochen}. Using (weak) toroidalization of the morphism $f \otimes \Q$, and induction on the dimension of $X\otimes \Q$ (to take care of the exceptional loci of the modifications used to obtain a toroidalization), one easily reduces to the following case. The morphism $f \otimes \Q$ is toroidal, $X \otimes \Q$ and  $Y \otimes \Q$ are nonsingular, and the zero loci and polar loci of $x_1, \dots, x_r$, restricted to $X \otimes \Q$, are contained in the support of the toroidal divisor on $X \otimes \Q$. Then $f \otimes \Q$ is log-smooth with respect to the toroidal divisors. In that case the Tameness Theorem follows directly from a logarithmic version of Hensel's lemma. We refer to \cite{Denef-AxKochen} for the details.
Note that the last sentence in the statement of the Tameness Theorem \ref{th-tameness} is a direct formal consequence of Lemma \ref{LemmaResidues1}.
The relation with logarithmic geometry is investigated in \cite{Cao-ColliotConjecture}.

The following Surjectivity Criterium is based on the Tameness Theorem \ref{th-tameness} and is essential for the proof of the Diophantine Excision Theorem \ref{th-purity}.

\begin{sectionItem} {\rm\textbf{Surjectivity Criterium.}} \label{surjectCriterium}
Let $f: X \rightarrow Y$ be a morphism of varieties over $\Z$, with $Y$ affine and $Y\otimes \Q$ smooth. Suppose that, given any regular rational functions $y_1, \dots, y_s$ on $Y$ and $M \in \N$, we have the following for almost all primes $p$. For each $b \in Y(\Z_p)$, with $\ord_p(y_i(b)) \leq M$ for $i = 1, \dots, s$, there exists $a \in X(\Z_p)$ such that $f(a)$ and $b$ have the same residues with respect to $y_1, \dots, y_s$.
\\ If this condition is satisfied, then the map $ X(\Z_p) \rightarrow Y(\Z_p)$, induced by $f$, is surjective for almost all primes $p$.
\end{sectionItem}

\noindent
\emph{Proof.} By the Tameness Theorem \ref{th-tameness} applied to the morphism $f$ and an empty list of rational functions on $X$, there exist a natural number $N_1$ and nonzero regular rational functions $y_1, \dots, y_s$ on $Y$ satisfying the conclusion of the Tameness Theorem for all primes $p > N_1$. By enlarging the list $y_1, \dots, y_s$ we may assume that it contains a set of affine coordinates for $Y$.  

Next we apply Lemma \ref{lem-forSurjectCrit} firstly to the regular rational functions  $y_1, \dots, y_s$ on $Y$, and secondly also to the regular rational functions  $y_1 \circ f, \dots, y_s \circ f$ on $X$. This yields a natural number $N_2$ and a common partition $\mathcal{P}$ of $\N^s$ satisfying, for all primes $p > N_2$, the conclusion of Lemma \ref{lem-forSurjectCrit} both for the $y_i$ on $Y$ and for the $y_i \circ f$ on $X$.

Choose a point in each stratum of this partition $\mathcal{P}$, and  choose $M \in \N$ bigger than the $p$-adic valuations of the coordinates of these points. Choose a natural number $N_3$ such that the hypothesis of the Surjectivity Criterium holds for the above $y_1, \dots, y_s$ and $M$, for all primes $p > N_3$. 

From now on, let $p$ be any prime bigger than $N_1, N_2$ and $N_3$, and let $b' \in Y(\Z_p)$. In order to prove the Surjectivity Criterium we will find an $a' \in X(\Z_p)$ with $f(a') = b'$. Because of observation \ref{obs-relatedToModifications}.(a), we may assume that $y_i(b') \ne 0$ for $i=1, \dots, s$. 
Set $z' := (y_1(b'), \dots, y_s(b'))$.

By our first application of Lemma \ref{lem-forSurjectCrit} and the choice of $M$, there exists a point $b \in Y(\Z_p)$ such that the $p$-adic valuations of $z := (y_1(b), \dots, y_s(b))$ and $z'$ are in the same stratum of the partition $\mathcal{P}$, and $\overline{\ac}(y_i(b)) = \overline{\ac}(y_i(b'))$, and $\ord_p(y_i(b)) \leq M$ for $i = 1, \dots, s$.
 
By the above mentioned instance of the hypothesis of the Surjectivity Criterium, there exists $a \in X(\Z_p)$ such that $f(a)$ and $b$ have the same residues with respect to $y_1, \dots, y_s$. 
Thus the $p$-adic valuations of $z'' := ( y_1(f(a)), \dots, y_s(f(a)) )$ and $z$ are equal and hence in the same stratum of $\mathcal{P}$ as these of $z'$.  
Moreover 
$\overline{\ac}(y_i(f(a))) = \overline{\ac}(y_i(b)) = \overline{\ac}(y_i(b'))$. 

Thus, by our second application of Lemma  \ref{lem-forSurjectCrit}, the multiplicative residue of $z'$ is realizable with respect to $y_1 \circ f, \dots, y_s \circ f$, because obviously the multiplicative residue of $z''$ is realizable with respect to these functions.  This means that there exists an $a'' \in X(\Z_p)$ such that $y_i(f(a''))$ and $y_i(b')$ have the same multiplicative residue for $i = 1, \dots, s$. Since the list $y_1, \dots, y_s$ contains a set of affine coordinates for $Y$, this implies that $f(a'') \equiv b'\!\mod p$. Hence $f(a'')$ and $b'$ have the same residues with respect to $y_1, \dots, y_s$. 

By our above mentioned application of the Tameness Theorem \ref{th-tameness} we conclude that there exists an $a' \in X(\Z_p)$ with $f(a') = b'$. This terminates the proof of the Surjectivity Criterium. 
$\square$

\section{The Diophantine Excision Theorem}
\label{sec-purity}

\begin{named} {Diophantine Excision Theorem} \label{th-purity} Let $f: X \rightarrow Y$ be a proper dominant morphism of varieties over $\Z$, with $Y\otimes \Q$ nonsingular. Assume that for each strict modification $f': X' \rightarrow Y'$ of $f$, with $Y'\otimes \Q$ nonsingular, there exists a closed subscheme $S'$ of $Y'$, of codimension $\geq 2$, such that for almost all primes $p$ we have
$$ \{y \in Y'(\Zp) \;|\; y \, \smod p \, \not \in  S'(\F_p) \} \subset f'(X'(\Zp)).$$
Then the map $X(\Zp) \rightarrow Y(\Zp)$, induced by $f$, is surjective for almost all primes $p$.
\end{named}
We prove the Diophantine Excision Theorem at the end of the present section, after two lemma's. But first we mention some observations whose proofs are straightforward.

\begin{sectionItem}\rm {\textbf{Observations}} \label{obs-relatedToExcision}

(a) Let $f: X \rightarrow Y$ be a proper dominant morphism of varieties over $\Z$, with $Y\otimes \Q$ nonsingular, satisfying the assumption of the Diophantine Excision Theorem. If $U$ is a nonempty open subscheme of $Y$, then also the morphism $ f^{-1}(U) \rightarrow U$, induced by $f$, satisfies the assumption of the Diophantine Excision Theorem.  \\
\emph{Proof of observation (a).} It suffices to show that any modification $\beta_0: U' \rightarrow U$ of $U$, with $U'\otimes \Q$ nonsingular, factors as an open immersion $j: U' \rightarrow Y'$ composed with a modification  $\beta: Y' \rightarrow Y$ of $Y$, with $Y'\otimes \Q$ nonsingular, and $j(U') = \beta^{-1}(U)$. To achieve this, let $\beta_1: U' \rightarrow Y$ be the composition of $\beta_0$ with the inclusion $U \subset Y$. Apply Nagata's compactification theorem (see e.g. \cite{Lutke-Compactification}) to factorize $\beta_1$ as an open immersion $j_2: U' \hookrightarrow Y''$ composed with a proper morphism $\beta_2: Y'' \rightarrow Y$ of $\Z$-varieties. This implies that $\beta_2$ is a modification of $Y$ and that $U'' := j_2(U') = \beta_2^{-1}(U)$. Indeed the composition of the morphisms 
$U' \hookrightarrow \beta_2^{-1}(U) \rightarrow U$, induced by $j_2$ and $\beta_2$, is proper since it equals $\beta_0$; thus the open immersion $U' \hookrightarrow \beta_2^{-1}(U)$ is proper and hence surjective. 
There exists a resolution of $Y'' \otimes \Q$ which is a composition of blowups with nonsingular centers $C_1, \dots, C_r$ that lie above $(Y'' \setminus U'') \otimes \Q$. Denote by $\bar{C}_1$ the closure of $C_1$ in $Y''$. Denote by $\bar{C}_2$ the closure of $C_2$ in the blowup of $Y''$ with center $\bar{C}_1$, and so on. Let $\pi: Y' \rightarrow Y''$ be the modification of $Y''$ obtained by composing the blowups with centers $\bar{C}_1, \dots, \bar{C}_r$. Then $Y' \otimes \Q$ is nonsingular. Moreover, $\pi$ is an isomorphism above $U''$. This yields an open immersion $j_1: U''\rightarrow Y'$. Clearly $j:= j_1 \circ j_2$ and $\beta := \beta_2 \circ \pi$ satisfy the required properties. $\square$

(b) Let $f: X \rightarrow Y$ be a proper dominant morphism of varieties over $\Z$, with $Y\otimes \Q$ nonsingular, satisfying the assumption of the Diophantine Excision Theorem. Let $f_1: X_1 \rightarrow Y_1$ be a strict modification of $f$, with $Y_1\otimes \Q$ nonsingular. Then $f_1$ satisfies the assumption of the Diophantine Excision Theorem. This follows directly from the fact that any strict modification of $f_1$ is also a strict modification of $f$.
\end{sectionItem}

\begin{lemma}\label{lem-1-forPurity} Let $f: X \rightarrow Y$ be a proper dominant morphism of varieties over $\Z$, with $Y$ smooth over $\Z$, satisfying the assumption of the Diophantine Excision Theorem.
Let $h: Y \rightarrow \A^1_{\Z}$ be a smooth morphism. Then for almost all primes $p$ we have the following: for each $b \in Y(\Z_p)$ there exists an $a \in X(\Z_p)$ such that $f(a) \equiv b \; \smod p$ and $h(f(a)) = h(b)$.
\end{lemma}
\noindent
\emph{Proof.} By noetherian induction, it suffices to show that for any integral closed subscheme $W$ of $Y$, there exists a nonempty open subscheme $W_0$ of $W$, such that, for almost all $p$, the assertion of the lemma holds for all $b \in Y(\Z_p)$ satisfying $b \smod p \in W_0$. Clearly, we may assume that $W \otimes \Q$ is nonempty, since otherwise $W(\F_p)$ is empty for almost all $p$. If $W = Y$ then we can directly apply the assumption of the Diophantine Excision Theorem, with $f' = f$, to find $W_0$. Hence we can assume that $W \varsubsetneq Y$.
Moreover, by cutting away the non-smooth locus of $W$ and using observation \ref{obs-relatedToExcision}.(a), we may also assume that $W$ is smooth over $\Z$.
Thus locally the ideal sheaf of W on Y can be generated by part of a set of local uniformizing parameters over \Z, hence its blowup is smooth over $\Z$ with exceptional locus a projective space bundle on $W$.

Let $\beta : Y' \rightarrow Y$ be the blowup of $Y$ with center $W$, and let $f': X' \rightarrow Y'$ be the strict transform of $f$ with respect to $\beta$. Because the assumption of the Diophantine Excision Theorem is assumed, there exists a closed subscheme $S'$ of $Y'$, of codimension $\geq 2$, such that for almost all primes $p$ we have
\begin{equation} \label{equation-1}
 \{y \in Y'(\Zp) \;|\; y \, \smod p \, \not \in  S'(\F_p) \} \subset f'(X'(\Zp)).
\end{equation}
To start, we take $W_0$ equal to $W$, but later on we will replace $W_0$ by a smaller nonempty open subscheme of $W$ if necessary.

When the restriction of $h$ to $W$ is dominant, then making $W_0$ smaller if necessary, we may suppose that the restriction of $h$ to $W_0$ is smooth. Whence the restriction of $h \circ \beta$ to $\beta^{-1}(W_0)$ is smooth, because the morphism from $\beta^{-1}(W_0)$ to $W_0$, induced by $\beta$, is smooth. This implies that $h \circ \beta$ is smooth at each point of $\beta^{-1}(W_0)$, by the smoothness criterium for morphisms of smooth schemes (Th\'{e}or\`{e}me 17.11.1 in \cite{EGA-4-4}).

When the restriction of $h$ to $W$ is not dominant, then $h(W\otimes \Q) = \{P\}$, with $P$ a closed point of $\A^1_{\Q}$. Let $[P]$ be the prime divisor on $\A^1_{\Q}$ consisting of this point P with multiplicity one. Because $h$ is smooth, the multiplicity of $\beta^{-1}(W)\otimes \Q$, in the divisor $((h \circ \beta)\otimes \Q)^*([P])$ on $Y'\otimes \Q$, equals $1$. Let $C'$ be the intersection of $\beta^{-1}(W)\otimes \Q$ with the union of the other irreducible components of this divisor. Clearly $C'$ has codimension $\geq 2$ in $Y'\otimes \Q$, and $(h \circ \beta)\otimes \Q$ is smooth at each point of $\beta^{-1}(W)\otimes \Q \setminus C'$. Enlarging $S'$ if necessary, we may suppose that $C' \subset S'$. Whence, the singular locus of $h \circ \beta$ is disjoint from $(\beta^{-1}(W) \setminus S')\otimes \Q$. Hence shrinking $W_0$ if necessary, by inverting finitely many primes, we can assume that the singular locus of $h \circ \beta$ is disjoint from $\beta^{-1}(W_0) \setminus S'$.

Thus in either case, we can suppose that $h \circ \beta$ is smooth at each point of $\beta^{-1}(W_0) \setminus S'$.

Since $S' \subset Y'$ has codimension $\geq 2$ and $\beta^{-1}(W) \subset Y'$ has codimension 1, the morphism from $\beta^{-1}(W) \setminus S'$ to $W$, induced by $\beta$, is dominant. Hence, replacing $W_0$ if necessary by a smaller open subscheme of $W$, we may suppose that
$ W_0 \subset \beta( \, Y' \setminus S' \, )$.

Let $p$ be a big enough prime, and consider any $b \in Y(\Z_p)$ satisfying $\bar{b} := b \,\smod p \in W_0$. Because the scheme-theoretic fibre of $\beta$ over $\bar{b}$ is isomorphic to a projective space over $\F_p$, and its intersection with $S'$ is contained in a hypersurface (of this fibre) with degree bounded independently of $b$, and because $p$ is big enough, there exists a $\F_p$-rational point $\bar{b'}$ on $Y' \setminus S'$ with $\beta (\bar{b'}) = \bar{b}$. Since $h \circ \beta$ is smooth at $\bar{b'}$ and $(h \circ \beta)(\bar{b'}) = h(b)  \smod p$, this point lifts to a point $b' \in Y'(\Z_p)$ with $b' \smod p = \bar{b'}$ and $(h \circ \beta)(b') = h(b)$. By (\ref{equation-1}) and because  
 $b' \smod p \not \in S'$, there exists $a' \in X'(\Z_p)$ with $f'(a') = b'$. Let $a \in X(\Z_p)$ be the image of $a'$ under the natural morphism $X' \rightarrow X$. Then $f(a) = 
 \beta (f'(a')) = \beta (b') \equiv b \; \smod p$, and $h(f(a)) = h(\beta (b')) = h(b)$. This terminates the proof of the lemma. $\square$

\begin{remark}\label{rem-1-forPurity} \rm The previous Lemma \ref{lem-1-forPurity} also holds when there is no $h$ involved, if we drop the requirement that $h(f(a)) = h(b)$. This follows formally from this lemma, using observation \ref{obs-relatedToExcision}.(a), by covering $Y$ by finitely many small enough open subschemes on which there exists a smooth morphism to $\A^1_{\Z}$.
\end{remark}

\begin{lemma}\label{lem-2-forPurity} Let $f: X \rightarrow Y$ be a proper dominant morphism of varieties over $\Z$, with $Y \otimes \Q$ nonsingular, satisfying the assumption of the Diophantine Excision Theorem.  Let $y_1, \dots, y_s$ be regular rational functions on $Y$, and $M \in \N$. Then for almost all primes $p$ we have the following. For each $b \in Y(\Z_p)$, with $\ord_p(y_i(b)) \leq M$ for $i = 1, \dots, s$, there exists $a \in X(\Z_p)$ such that $f(a)$ and $b$ have the same residues with respect to $y_1, \dots, y_s$.
\end{lemma}
\noindent
\emph{Proof.} 
Let $D$ be the union of the zero loci of the regular rational functions $y_1, \dots, y_s$ on $Y$, considered as a subset of $Y$. Using embedded resolution of singularities of $D  \otimes \Q \subset Y \otimes \Q$, modifying $Y$, without changing  $Y \otimes \Q \setminus D \otimes \Q$, replacing $f$ by its strict transform with respect to the modification of $Y$, and inverting a finite number of primes, we may assume the following. The variety $Y$ is smooth over $\Z$, and $D$ is a reduced strict normal crossings divisor over $\Z$ (in the sense of subsection \ref{conventions}). This reduction is easily verified applying the valuative criterion of
properness to the resolution morphism and using observation \ref{obs-relatedToExcision}.(b). Thus, by covering $Y$ with finitely many suitable open subschemes, and using observation \ref{obs-relatedToExcision}.(a), we can further assume that $Y$ is affine, and that each $y_i$ can be written as a unit in $\Gamma(Y, \mathcal{O}_Y)$ times a monomial in uniformizing parameters over $\Z$ on $X$ (i.e. regular functions on $Y$ that induce an \'{e}tale morphism to an affine space over $\Z$, see subsection \ref{conventions}).
Hence, by Lemma \ref{LemmaResidues2}, we can moreover assume that $y_1, \dots, y_s$ are part of a set of uniformizing parameters over $\Z$ on $Y$.

It remains now to prove the lemma in the special case that $Y$ is smooth over $\Z$, and affine, say $Y = \Spec (A)$, and that $y_1, \dots, y_s$ are part of a set of uniformizing parameters over $\Z$ on $Y$. We prove this special case by induction on $M$. Let $p$ be a prime, big enough with respect to $M$ and all data, and let $b \in Y(\Z_p)$ be any point with $\ord_p(y_i(b)) \leq M$ for all $i = 1, \dots, s$. If $M = 0$, then, in order to prove the lemma, it suffices to find $a \in X(\Z_p)$ with $f(a) \equiv b \; \smod p$. The existence of such an $a$ follows from Remark \ref{rem-1-forPurity}.
Thus we may suppose that $M > 0$ and that
$$ I_0 := \{i \in \N \; | \; \ord_p (y_i(b)) > 0\, , \, 1 \leq i \leq s \} \not = \emptyset \, . $$
Choose $i_0 \in I_0$ such that $\ord_p (y_{i_0}(b)) = \Min_{i \in I_0} \ord_p (y_i(b))$.

Let $\pi: Y' \rightarrow Y$ be the blowup of the ideal sheaf on $Y$ generated by all the $y_i$ with $i \in I_0$. Consider the chart $U$ on $Y'$, defined as follows:
$$ U:= \Spec(A[(y_i/y_{i_0})_{i \in I_0}]) \, \xrightarrow{\pi} \, \Spec(A) = Y \, . $$
There exists a unique $b' \in U(\Z_p)$ with $\pi(b') = b$. Set $y'_i = y_i/y_{i_0}$ for $i \in I_0 \setminus \{ i_0 \}$, and $y'_i = y_i$ for the other $i \in \{1, \dots, s \}$. One easily verifies that $y'_1, \dots, y'_s$ are part of a set of uniformizing parameters over $\Z$ on $U$. Clearly, either $0 \leq \ord_p(y'_i(b')) < M$, for all $i$, or $\ord_p(y'_i(b')) = 0$, for all $i \not = i_0$. We call these respectively the first case and the second case.

Let $f': X' \rightarrow Y'$ be the strict transform of $f$ with respect to the blowup $\pi: Y' \rightarrow Y$. In the first case, we apply the induction hypothesis to the morphism $f'^{-1}(U) \rightarrow U$ induced by $f'$, and the regular rational functions $y'_1, \dots, y'_s$ on $U$, to find $a' \in X'(\Z_p)$, with $f'(a') \in U(\Z_p)$, such that $f'(a')$ and $b'$ have the same residues with respect to these functions on $U$. In the second case, we apply Lemma \ref{lem-1-forPurity} to the morphism $f'^{-1}(U) \rightarrow U$ induced by $f'$, and the morphism $U \rightarrow \A^1_{\Z}$ induced by $y_{i_0}$, to find $a' \in X'(\Z_p)$, with $f'(a') \in U(\Z_p)$, such that $f'(a') \equiv b' \; \smod p$ and $y_{i_0}(f'(a')) = y_{i_0}(b')$. Hence, also in the second case, $f'(a')$ and $b'$ have the same residues with respect to $y'_1, \dots, y'_s$, because $y'_i(b')$ is a unit in $\Z_p$ for all $i \not = i_0$.

Denote by $a$ the image of $a'$ under the natural map $X'(\Z_p) \rightarrow X(\Z_p)$. Then the points $f(a) = \pi(f'(a'))$ and $b = \pi(b')$ have the same residues with respect to $y_1, \dots, y_s$. This terminates the proof of the lemma. $\square$

\begin{sectionItem} Proof of the Diophantine Excision Theorem \ref{th-purity}
\end{sectionItem} \noindent
Using observation \ref{obs-relatedToExcision}.(a) we may assume that $Y$ is affine. Then the Excision Theorem \ref{th-purity} is a direct consequence of the above Lemma \ref{lem-2-forPurity} and the Surjectivity Criterium \ref{surjectCriterium}. $\square$


\section{Proof of the Main Theorem \ref{th-main}}  \label{sec-proofMainThm}

\noindent
In this section we show that the Main Theorem \ref{th-main} is an easy consequence of the Diophantine Excision Theorem \ref{th-purity} and the following lemma whose proof is rather straightforward.

\begin{lemma}\label{lem-1-forMainThm} Let $f: X \rightarrow Y$ be a proper dominant morphism of smooth varieties over $\Z$, with geometrically integral generic fibre. Assume for each $\Z$-flat prime divisor $D$ on $Y$, that the divisor $f^{*}(D)$ on $X$ has an irreducible component $C$ with multiplicity 1 and geometrically integral generic fibre over $D$ (i.e. the morphism $C \rightarrow D$, induced by $f$, has geometrically integral generic fibre).
Then there exists a closed subscheme $S$ of $Y$, of codimension $\geq 2$, such that for almost all primes $p$ we have
$$ \{y \in Y(\Zp) \;|\; y \, \smod p \, \not \in  S(\F_p) \} \subset f(X(\Zp)).$$
\end{lemma}
\noindent
\textit{Proof.} By Th\'{e}or\`{e}me 9.7.7 of \cite{EGA-4-3} , there exists a reduced closed subscheme $E \subset Y$, of pure codimension 1, such that over the complement of $E$, the morphism $f$ is smooth with geometrically integral fibres. Hence, for almost all primes $p$, any $y \in Y(\Zp)$, with $y \, \smod p \, \not \in  E(\F_p)$, belongs to $f(X(\Zp))$. Indeed this follows from Hensel's Lemma and the Lang-Weil bound \cite{Lang-Weil}.

For each irreducible component $D$ of $E$ we reason as follows. If $D$ is not flat over $\Spec(\Z)$, then $D(\F_p)$ is empty for almost all primes $p$. Suppose now that $D$ is flat over $\Spec(\Z)$. By assumption, the divisor $f^{*}(D)$ on $X$ has an irreducible component $C$ with multiplicity 1 and geometrically integral generic fibre over $D$. In particular, $C$ dominates $D$. Hence there exists a reduced closed subscheme $S$ of $D$, of codimension $\geq 1$ in $D$, such that, over the complement of $S$, all fibres of $C \xrightarrow{f} D$ are geometrically integral and intersect the smooth locus of $f: X \rightarrow Y$. Indeed, $f$ is smooth at the generic point of $C$, because $C$ has multiplicity 1 in the divisor $f^{*}(D)$. Again by Hensel's Lemma and the Lang-Weil bound \cite{Lang-Weil}, we conclude for almost all primes $p$ that any $y \in Y(\Zp)$, with $y \, \smod p \,  \in  D(\F_p) \setminus S(\F_p)$, belongs to $f(X(\Zp))$.

Taking the union of the subschemes $S$, obtained as above for each $\Z$-flat irreducible component $D$ of $E$, we obtain a closed subscheme of $Y$, of codimension $\geq 2$, that satisfies the conclusion of the lemma. $\square$

\begin{sectionItem} Proof of the Main Theorem \ref{th-main}
\end{sectionItem}
Let $f: X \rightarrow Y$ be a dominant morphism of nonsingular proper geometrically integral varieties over $\Q$, which satisfies the hypotheses of the Main Theorem. Choose a proper dominant morphism $\tilde{f}: \tilde{X} \rightarrow \tilde{Y}$ of smooth varieties over $\Z$, whose base change to $\Q$ is isomorphic to $f$. Because $X$ and $Y$ are proper, it suffices to prove that the map $X(\Z_p) \rightarrow Y(\Z_p)$, induced by $\tilde{f}$, is surjective for almost all primes $p$. By the Diophantine Excision Theorem \ref{th-purity}, it suffices to prove that the morphism $\tilde{f}$ satisfies the assumption in the Diophantine Excision Theorem, with $f$ replaced by $\tilde{f}$.

Let $\tilde{f}': \tilde{X}' \rightarrow \tilde{Y}'$ be any strict modification of $\tilde{f}$, with $\tilde{Y}'\otimes \Q$ nonsingular. Note that the generic fibre of  $\tilde{f}'$ equals the one of $\tilde{f}$ and is contained in the smooth locus of $\tilde{X}'$, because the modification is strict and $\tilde{X}$ is smooth.
We have to prove that there exists a closed subscheme $S$ of $\tilde{Y}'$, of codimension $\geq 2$, such that for almost all primes $p$ we have
$$ \{y \in \tilde{Y}'(\Zp) \;|\; y \, \smod p \, \not \in  S(\F_p) \} \subset \tilde{f}'(\tilde{X}'(\Zp)).$$
Composing $\tilde{f}'$ with a morphism whose base change to $\Q$ resolves the singularities of $\tilde{X}'\otimes \Q$ without changing the smooth locus of $\tilde{X}'\otimes \Q$, and inverting a finite number of primes, we see that in order to prove the above, we may assume  the following. The varieties $\tilde{X}'$ and $\tilde{Y}'$ are smooth over $\Z$, and $\tilde{f}'$ is a modification of $\tilde{f}$, with the same generic fibre as $\tilde{f}$. But now $\tilde{f}'$ is not necessarily a strict modification of $\tilde{f}$ anymore.

Because, by assumption, $f$ satisfies the hypotheses of the Main Theorem, it is easy to verify that $\tilde{f}'$ satisfies the hypotheses of Lemma \ref{lem-1-forMainThm} (with $f$ replaced by $\tilde{f}'$). Hence this lemma implies the existence of a closed subscheme $S$ of $\tilde{Y}'$ with the required properties. This terminates the proof of the Main Theorem.  $\square$

\begin{sectionItem} An alternative proof of Colliot-Th\'{e}l\`{e}ne's Conjecture. \label{prf-alternative}
\end{sectionItem}
Using model theory (mathematical logic) one can give a much simpler proof of Colliot-Th\'{e}l\`{e}ne's Conjecture \ref{conj-Colliot}. However we don't see how to extend this to prove the stronger Theorem \ref{th-main} or the Diophantine Excision Theorem \ref{th-purity}. Moreover one of the motivations of Colliot-Th\'{e}l\`{e}ne was to obtain a new proof of the Ax-Kochen Theorem which does not rely on methods from mathematical logic. We briefly sketch this simpler proof of Colliot-Th\'{e}l\`{e}ne's Conjecture.

Assume the notation and hypotheses in the formulation of Conjecture \ref{conj-Colliot}. Using a standard argument from model theory and the Ax-Kochen-Er\v sov Transfer Principle \cite{AxKoch-1, Ersov-1}, we will first show that, in order to prove the conjecture, it suffices to show that the map $X(F[[t]]) \rightarrow Y(F[[t]])$, induced by $f$, is surjective for any pseudo algebraically closed field $F$ of characteristic zero. This goes as follows. If Conjecture \ref{conj-Colliot} is false for $f$, then there exists an infinite set $S$ of primes $p$ for which the map from $X(\Q_p)$ to $Y(\Q_p)$, induced by $f$, is not surjective. Let $K$ be the ultra product of all the fields $\Q_p$ with respect to an ultra filter, on the set $\mathcal{P}$ of all primes, containing $S$ and each subset of $\mathcal{P}$ with finite complement. Then the map from $X(K)$ to $Y(K)$, induced by $f$, is not surjective. Notice that $K$ is a henselian valued field with residue field a pseudo algebraically closed field $F$ of characteristic zero (by \cite{Lang-Weil}), and value group elementary equivalent to $\Z$. Hence $K$ is elementary equivalent to the field of fractions $F((t))$ of $F[[t]]$, by the Ax-Kochen-Er\v sov Transfer Principle which states that any two henselian valued fields are elementary equivalent if they have elementary equivalent value groups and elementary equivalent residue fields of characteristic zero. Thus the map from $X(F((t)))$ to $Y(F((t)))$, induced by $f$, is not surjective if Conjecture \ref{conj-Colliot} would be false for $f$. Since $f$ is proper the same holds for $F((t))$ replaced by $F[[t]]$. We conclude that in order to prove the conjecture, it suffices to show that the map $X(F[[t]]) \rightarrow Y(F[[t]])$, induced by $f$, is surjective for any pseudo algebraically closed field $F$ of characteristic zero.

Let $y \in Y(F[[t]])$. We have to show that $y \in f(X(F[[t]]))$. Let $s$ be the closed point of $\Spec(F[[t]])$. By slightly moving $y$ and using Greenberg's Theorem \cite{Greenberg}, we may assume that the homomorphism $\cO_{Y,y(s)} \rightarrow F[[t]]$ induced by $y$ is injective. Composing this homomorphism with the standard valuation on $F[[t]]$, induces a discrete valuation $\nu$ on the function field $K$ of $Y$, with valuation ring say $A$.

If $\nu$ is trivial, then $y(s)$ is the generic point of $Y$. Hence $f$ is smooth at each point in the fibre of $y(s)$, and this fibre is geometrically integral. This implies that $y$ lifts to a $F[[t]]$-rational point $x$ on $X$, by Hensel's Lemma and the assumption that $F$ is pseudo algebraically closed.

Thus we may assume that the discrete valuation $\nu$ is not trivial. Hence there exists an integral regular $A$-scheme $\fX$ as in the formulation of Conjecture \ref{conj-Colliot}. Note that $y$ induces a $F[[t]]$-rational point $\tilde{y}$ on $\Spec(A)$, and a homomorphism $K \rightarrow F((t))$.  Using the hypothesis about the special fibre of $\fX$, Hensel's Lemma, and the assumption that $F$ is pseudo algebraically closed, one easily verifies that $\tilde{y}$ lifts to a $F[[t]]$-rational point $\tilde{x}$ on $\fX$. Because the generic fibre of $\fX$ is $K$-isomorphic to the generic fibre of $f$, and because $\tilde{x}$ extends to a $F((t))$-rational point on $\fX \otimes K$, we find a $F((t))$-rational point on $X$, and hence, by the properness of $X$, also a $F[[t]]$-rational point $x$ on $X$ with $f(x) = y$.  $\square$

\bibliographystyle{plain}
\bibliography{JD}

\end{document}